\documentclass[10pt]{article}
\usepackage{cite}
\usepackage{mathrsfs}
\usepackage{amsfonts}
\usepackage{amsmath}
\usepackage{amsfonts,amssymb,color}
\usepackage{dsfont}
\usepackage{curves}
\usepackage{mathrsfs}
\usepackage{pifont}
\usepackage{amssymb}
\allowdisplaybreaks

\numberwithin{equation}{section}

\date{}

\def\BigRoman{\uppercase\expandafter{\romannumeral\number\count 255 }}
\def\Romannumeral{\afterassignment\BigRoman\count255=}

\begin{document}
\title{Generalized toughness and $Q$-index in a graph
}
\author{\small Sizhong Zhou\footnote{Corresponding author. E-mail address: zsz\_cumt@163.com (S. Zhou)}\\
\small  School of Science, Jiangsu University of Science and Technology,\\
\small  Zhenjiang, Jiangsu 212100, China\\
}

\maketitle
\begin{abstract}
\noindent Let $G$ be a graph. We denote by $c(G)$, $\alpha(G)$ and $q(G)$ the number of components, the independence number and the signless Laplacian spectral radius ($Q$-index for short) of $G$,
respectively. The toughness of $G$ is defined by $t(G)=\min\left\{\frac{|S|}{c(G-S)}:S\subseteq V(G), c(G-S)\geq2\right\}$ for $G\neq K_n$ and $t(G)=+\infty$ for $G=K_n$. Chen, Gu and Lin [Generalized
toughness and spectral radius of graphs, Discrete Math. 349 (2026) 114776] generalized this notion and defined the $l$-toughness $t_l(G)$ of a graph $G$ as $t_l(G)=\min\left\{\frac{|S|}{c(G-S)}:S\subset
V(G), c(G-S)\geq l\right\}$ if $2\leq l\leq\alpha(G)$, and $t_l(G)=+\infty$ if $l>\alpha(G)$. If $t_l(G)\geq t$, then $G$ is said to be $(t,l)$-tough.
In this paper, we put forward $Q$-index conditions for a graph to be $(b,l)$-tough and $(\frac{1}{b},l)$-tough, respectively.\\
\begin{flushleft}

(2020) Mathematics Subject Classification: 05C50, 05C35
\end{flushleft}
\end{abstract}

\section{Introduction}

Throughout this paper, we deal only with undirected and simple graphs. For any graph $G$, its vertex set and edge set are denoted by $V(G)$ and $E(G)$, respectively. We write $|V(G)|=n$ and $|E(G)|=e(G)$
for the order and the size of $G$, respectively. Let $i(G)$ and $\alpha(G)$ denote the number of isolated vertices and the independence number of $G$, respectively. For any $v\in V(G)$, let $d_G(v)$ denote
the degree of $v$ in $G$. Let $S$ be a subset of $V(G)$. We write $G[S]$ and $G-S$ for the subgraphs of $G$ induced by $S$ and $V(G)\setminus S$, respectively. The cycle and the complete graph of order $n$
are denoted by $C_n$ and $K_n$, respectively. Given two graphs $G_1$ and $G_2$, we denote by $G_1\cup G_2$ their disjoint union and by $tG_1$ the disjoint union of $t$ copies of $G_1$. The join $G_1\vee G_2$
is the graph obtained from $G_1\cup G_2$ by adding all possible edges between $V(G_1)$ and $V(G_2)$.

For a given graph $G$ with $V(G)=\{v_1,v_2,\ldots,v_n\}$, the adjacency matrix of $G$ is the 0-1 matrix $A(G)=(a_{ij})_{n\times n}$, where $a_{ij}=1$ if the two vertices $v_i$ and $v_j$ are adjacent and
$a_{ij}=0$ otherwise. The spectral radius of $G$, denoted by $\rho(G)$, is the largest eigenvalue of $A(G)$. We define $Q(G)=D(G)+A(G)$ for the signless Laplacian matrix of $G$, where
$D(G)=\mbox{diag}(d_G(v_1),d_G(v_2),\ldots,d_G(v_n))$ is the diagonal degree matrix of $G$. The signless Laplacian spectral radius ($Q$-index for short) of $G$, denoted by $q(G)$, is the largest eigenvalue
of $Q(G)$.

In the past decades, lots of scholars established the connections between $Q$-index and the structural properties of graphs. Pan, Li and Zhao \cite{PLZ} established the relationships between the fractional
matching number and $Q$-index of a graph, and provided some sufficient conditions with respect to $Q$-index for a graph to possess a fractional perfect matching. Hao and Li \cite{HL} proposed some lower
bounds based on $Q$-index to guarantee that a graph contains a path-factor. Jahanbani and Sheikholeslami \cite{JS} presented the relationship between the sum-connectivity index and $Q$-index in a graph.
Zhang and Wang \cite{ZW} characterized the unique extremal graph with maximum $Q$-index among all $2K_3$-free graphs with order at least 44. Cvetkovi\'c and Simi\'c \cite{CS} obtained a spectral theory of
graphs based on $Q$-index and compared it with other spectral theories. Zhou and Zhang \cite{ZZ} determined a lower bound on $Q$-index of a graph $G$ to guarantee that $G$ is $k$-extendable and claimed
some extremal graphs to show that all the bounds are sharp. Zheng, Li, Luo and Wang \cite{ZLLW} created a sufficient condition by virtue of $Q$-index for a graph with minimum degree to be $k$-factor-critical.
Wu, Zhou and Liu \cite{WZL} provided a lower bound on $Q$-index of a bipartite graph $G$ with bipartition $(A,B)$, in which the bound ensures that $G$ contains a spanning tree $T$ with $d_T(v)\geq k$ for
all $v\in A$. Wu \cite{Ws} obtained a sufficient condition in view of $Q$-index to ensure that a graph has the $k$-strong parity property. Wu \cite{Wu} established a tight sufficient condition involving
$Q$-index to guarantee the existence of $\{C_{2i+1},T:1\leq i<\frac{r}{k-r},T\in\mathcal{T}_{\frac{k}{r}}\}$-factor in a graph, where $r<k$ are two positive integers and $\mathcal{T}_{\frac{k}{r}}$ denotes
the set of trees $T$ such that $i(T-S)\leq\frac{k}{r}|S|$ for any $S\subset V(T)$ and for any $e\in E(T)$ there exists a set $S^{*}\subset V(T)$ with $i((T-e)-S^{*})>\frac{k}{r}|S^{*}|$. For the results
on other spectral radius in a graph, we refer the readers to \cite{FS,DM,MK,Zs,ZW1,ZZL,ZZL1,Wc,WZ,ZL}.

In 1973, Chv\'atal \cite{Ch} introduced the notion of toughness $t(G)$: If $G$ is complete, $t(G)=+\infty$. If $G$ is not complete,
$$
t(G)=\min\left\{\frac{|S|}{c(G-S)}:S\subset V(G), c(G-S)\geq2\right\},
$$
where $c(G-S)$ denotes the number of components of $G-S$. If $t(G)\geq t$, then $G$ is said to be $t$-tough. The toughness is related to many properties of a graph. Chv\'atal \cite{Ch} created some
connections between toughness and the existence of Hamiltonian cycles. Enomoto, Jackson, Katerinis and Saito \cite{EJKS} provided a toughness condition for a graph to possess a $k$-factor. Liu and Zhang
\cite{LZ} showed a toughness condition for the existence of a fractional $k$-factor in a graph. Fan, Lin and Lu \cite{FLL} gave a spectral sufficient condition for a connected graph with minimum degree
$\delta\geq2$ to be 1-tough. Jia and Lou \cite{JL} investigated the relationships between $Q$-index and tough graphs. Zhou \cite{Zt} posed distance spectral radius conditions for a connected graph to be
$t$-tough and $\frac{1}{t}$-tough, respectively.

Chen, Gu and Lin \cite{CGL} generalized the toughness notion and defined the $l$-toughness $t_l(G)$ of a graph $G$ as follows:
$$
t_l(G)=\min\left\{\frac{|S|}{c(G-S)}:S\subset V(G), c(G-S)\geq l\right\},
$$
where $2\leq l\leq\alpha(G)$. If $l>\alpha(G)$, then we may define $t_l(G)=+\infty$ by convention. According to definition, $t_2(G)$ is the classic toughness $t(G)$. If $t_l(G)\geq t$, then $G$ is said to
be $(t,l)$-tough. Furthermore, they arose spectral conditions for a graph to be $(b,l)$-tough and $(\frac{1}{b},l)$-tough, respectively.

Motivated by \cite{CGL} directly, we provide a lower bound on $Q$-index in a graph to guarantee that a graph is $(b,l)$-tough and $(\frac{1}{b},l)$-tough, respectively.

\medskip

\noindent{\textbf{Theorem 1.1.}} Let $b\geq1$ and $l\geq2$ be two integers, and let $G$ be a connected graph of order $n$ with $n\geq\max\{(\frac{5}{2}b^{2}+4b+3)l-b^{2}-2b-5,\frac{(2b+1)l^{2}+(2b-3)l+2}{2}\}$. If
$$
q(G)\geq q(K_{bl-1}\vee(K_{n-(b+1)l+2}\cup(l-1)K_1)),
$$
then $t_{l}(G)\geq b$, unless $G=K_{bl-1}\vee(K_{n-(b+1)l+2}\cup(l-1)K_1)$.

\medskip

\noindent{\textbf{Theorem 1.2.}} Let $b\geq2$ and $l\geq2$ be two integers, and let $G$ be a connected graph of order $n$ with $n\geq6b\left\lceil\frac{l-1}{b}\right\rceil$. If
$$
q(G)\geq q(K_{\lfloor\frac{l-1}{b}\rfloor}\vee(K_{n-\lfloor\frac{l-1}{b}\rfloor-l+1}\cup(l-1)K_1)),
$$
then $t_{l}(G)\geq\frac{1}{b}$, unless $G=K_{\lfloor\frac{l-1}{b}\rfloor}\vee(K_{n-\lfloor\frac{l-1}{b}\rfloor-l+1}\cup(l-1)K_1)$.

\section{Some preliminaries}

Let $M$ be a real matrix of order $n$ and $\mathcal{N}=\{1,2,\ldots,n\}$. Given a partition $\pi:\mathcal{N}=\mathcal{N}_1\cup\mathcal{N}_2\cup\cdots\cup\mathcal{N}_r$, the matrix $M$ can
be correspondingly described as
\begin{align*}
M=\left(
  \begin{array}{cccc}
    M_{11} & M_{12} & \cdots & M_{1r}\\
    M_{21} & M_{22} & \cdots & M_{2r}\\
    \vdots & \vdots & \ddots & \vdots\\
    M_{r1} & M_{r2} & \cdots & M_{rr}\\
  \end{array}
\right).
\end{align*}
The quotient matrix of $M$ with respect to the partition $\pi$ is defined by the matrix $M_{\pi}=(m_{ij})$ of order $r$, where $m_{ij}$ denotes the average row sum of the block $M_{ij}$. If for any $i,j$, the block
$M_{ij}$ of $M$ has constant row sum, then the partition $\pi$ is called equitable.

\medskip

\noindent{\textbf{Lemma 2.1}} (You, Yang, So and Xi \cite{YYSX}). Let $M$ be a real matrix of order $n$ with an equitable partition $\pi$, and let $M_{\pi}$ be the corresponding quotient matrix. Then the eigenvalues
of $M_{\pi}$ are also eigenvalues of $M$. Furthermore, if $M$ is nonnegative and irreducible, then the largest eigenvalues of $M$ and $M_{\pi}$ are equal.

\medskip

\noindent{\textbf{Lemma 2.2}} (Shen, You, Zhang and Li \cite{SYZL}). If $H$ is a subgraph of a connected graph $G$, then
$$
q(H)\leq q(G),
$$
where the equality occurs if and only if $H=G$.

\medskip

\noindent{\textbf{Lemma 2.3}} (Zheng, Li, Luo and Wang \cite{ZLLW}). Let $\sum\limits_{i=1}^{t}n_i=n-s$. If $n_1\geq n_2\geq\cdots\geq n_t\geq p\geq1$ and $n_1<n-s-p(t-1)$, then
$$
q(K_s\vee(K_{n_1}\cup K_{n_2}\cup\cdots\cup K_{n_t}))<q(K_s\vee(K_{n-s-p(t-1)}\cup(t-1)K_p)).
$$

\medskip

\noindent{\textbf{Lemma 2.4}} (Das \cite{Das}, Feng and Yu \cite{FY}). Let $G$ be a graph of order $n$. Then
$$
q(G)\leq\frac{2e(G)}{n-1}+n-2.
$$

\medskip

\section{The proof of Theorem 1.1}

\noindent{\it Proof of Theorem 1.1.} Suppose that $t_l(G)<b$ with $l\geq2$. Then we possess $|S|<bc(G-S)$ for some nonempty subset $S\subseteq V(G)$. Let $|S|=s$ and $c(G-S)=\omega$. Then $s\leq b\omega-1$.
The following proof will be divided into two cases in terms of the value of $n$.

\noindent{\bf Case 1.} $n\geq(b+1)\omega-1$.

In this case, we easily see that $G$ is a spanning subgraph of $G_1=K_{b\omega-1}\vee(K_{n_1}\cup K_{n_2}\cup\cdots\cup K_{n_{\omega}})$ for some positive integers $n_1\geq n_2\geq\cdots\geq n_{\omega}$ with
$\sum\limits_{i=1}^{\omega}n_i=n-b\omega+1$. Together with Lemma 2.2, we conclude
\begin{align}\label{eq:3.1}
q(G)\leq q(G_1),
\end{align}
where the equality occurs if and only if $G=G_1$. Write $G_2=K_{b\omega-1}\vee(K_{n-(b+1)\omega+2}\cup(\omega-1)K_1)$. According to Lemma 2.3, we get
\begin{align}\label{eq:3.2}
q(G_1)\leq q(G_2),
\end{align}
where the equality occurs if and only if $(n_1,n_2,\ldots,n_{\omega})=(n-(b+1)\omega+2,1,\ldots,1)$. For $\omega=l$, we possess $G_2=K_{bl-1}\vee(K_{n-(b+1)l+2}\cup(l-1)K_1)$. Together with \eqref{eq:3.1} and
\eqref{eq:3.2}, we obtain
$$
q(G)\leq q(K_{bl-1}\vee(K_{n-(b+1)l+2}\cup(l-1)K_1)),
$$
where the equality holds if and only if $G=K_{bl-1}\vee(K_{n-(b+1)l+2}\cup(l-1)K_1)$, which contradicts the condition of the theorem because $K_{bl-1}\vee(K_{n-(b+1)l+2}\cup(l-1)K_1)$ is not $(b,l)$-tough.
Next, we consider $\omega\geq l+1$.

Note that $K_{n-l+1}$ is a proper subgraph of $K_{bl-1}\vee(K_{n-(b+1)l+2}\cup(l-1)K_1)$. Combining this with Lemma 2.2, we deduce
\begin{align}\label{eq:3.3}
q(K_{bl-1}\vee(K_{n-(b+1)l+2}\cup(l-1)K_1))>q(K_{n-l+1})=2n-2l.
\end{align}

Recall that $G_2=K_{b\omega-1}\vee(K_{n-(b+1)\omega+2}\cup(\omega-1)K_1)$. By virtue of Lemma 2.4, we obtain
\begin{align}\label{eq:3.4}
q(G_2)\leq&\frac{2e(G_2)}{n-1}+n-2\nonumber\\
=&\frac{2\binom{n-\omega+1}{2}+2(b\omega-1)(\omega-1)}{n-1}+n-2\nonumber\\
=&\frac{(2b+1)\omega^{2}-(2n+2b+3)\omega+2n^{2}-2n+4}{n-1}.
\end{align}
Let $\varphi(\omega)=(2b+1)\omega^{2}-(2n+2b+3)\omega+2n^{2}-2n+4$. Recall that $n\geq(b+1)\omega-1$. Then we deduce $l+1\leq\omega\leq\frac{n+1}{b+1}$. By a simple computation, we possess
$$
\varphi(l+1)-\varphi\left(\frac{n+1}{b+1}\right)=\frac{1}{(b+1)^{2}}(n-(b+1)l-b)(n-(b+1)(2b+1)l+1)\geq0
$$
by $n\geq\max\{(\frac{5}{2}b^{2}+4b+3)l-b^{2}-2b-5,\frac{(2b+1)l^{2}+(2b-3)l+2}{2}\}\geq(\frac{5}{2}b^{2}+4b+3)l-b^{2}-2b-5\geq(b+1)(2b+1)l-1>(b+1)l+b$, which implies
\begin{align}\label{eq:3.5}
\varphi(\omega)\leq\varphi(l+1)
\end{align}
for $l+1\leq\omega\leq\frac{n+1}{b+1}$. It follows from \eqref{eq:3.1}--\eqref{eq:3.5} and $n\geq\max\{(\frac{5}{2}b^{2}+4b+3)l-b^{2}-2b-5,\frac{(2b+1)l^{2}+(2b-3)l+2}{2}\}\geq\frac{(2b+1)l^{2}+(2b-3)l+2}{2}$
that
\begin{align*}
q(G)\leq&q(G_1)\leq q(G_2)\leq\frac{\varphi(\omega)}{n-1}\leq\frac{\varphi(l+1)}{n-1}\\
=&\frac{(2b+1)(l+1)^{2}-(2n+2b+3)(l+1)+2n^{2}-2n+4}{n-1}\\
=&2n-2l-\frac{2n-(2b+1)l^{2}-(2b-3)l-2}{n-1}\\
\leq&2n-2l\\
<&q(K_{bl-1}\vee(K_{n-(b+1)l+2}\cup(l-1)K_1)),
\end{align*}
which contradicts $q(G)\geq q(K_{bl-1}\vee(K_{n-(b+1)l+2}\cup(l-1)K_1))$.

\noindent{\bf Case 2.} $n\leq(b+1)\omega-2$.

In this case, we easily see that $G$ is a spanning subgraph of graph $K_{n-\omega}\vee\omega K_1$ with $\omega\geq\lceil\frac{n+2}{b+1}\rceil$. Write $G_3=K_{n-\lceil\frac{n+2}{b+1}\rceil}\vee\lceil\frac{n+2}{b+1}\rceil K_1$.
Obviously, graph $K_{n-\omega}\vee\omega K_1$ is a spanning subgraph of $G_3=K_{n-\lceil\frac{n+2}{b+1}\rceil}\vee\lceil\frac{n+2}{b+1}\rceil K_1$. Combining these with Lemma 2.2, we claim
\begin{align}\label{eq:3.6}
q(G)\leq q(K_{n-\omega}\vee\omega K_1)\leq q(G_3),
\end{align}
with equalities holding if and only if $G=G_3$. Notice that $2e(G_3)=2\binom{n-\lceil\frac{n+2}{b+1}\rceil}{2}+2\lceil\frac{n+2}{b+1}\rceil(n-\lceil\frac{n+2}{b+1}\rceil)
=(n-\lceil\frac{n+2}{b+1}\rceil)(n+\lceil\frac{n+2}{b+1}\rceil-1)<(n-\frac{n+2}{b+1})(n+\frac{n+2}{b+1})=\frac{1}{(b+1)^{2}}((b^{2}+2b)n^{2}-4n-4)$. Combining this with Lemma 2.4, $l\geq2$ and $n\geq\max\{(\frac{5}{2}b^{2}+4b+3)l-b^{2}-2b-5,\frac{(2b+1)l^{2}+(2b-3)l+2}{2}\}\geq(\frac{5}{2}b^{2}+4b+3)l-b^{2}-2b-5>(2b^{2}+4b+2)l-b^{2}-2b-5$, we possess
\begin{align}\label{eq:3.7}
q(G_3)\leq&\frac{2e(G_3)}{n-1}+n-2\nonumber\\
<&\frac{(b^{2}+2b)n^{2}-4n-4}{(b+1)^{2}(n-1)}+n-2\nonumber\\
=&2n-2l-\frac{n^{2}-(2l(b^{2}+2b+1)-b^{2}-2b-5)n+2l(b^{2}+2b+1)-2b^{2}-4b+2}{(b+1)^{2}(n-1)}\nonumber\\
=&2n-2l-\frac{(n-1)(n-2l(b^{2}+2b+1)+b^{2}+2b+6)-b^{2}-2b+8}{(b+1)^{2}(n-1)}\nonumber\\
<&2n-2l.
\end{align}

According to \eqref{eq:3.3}, \eqref{eq:3.6} and \eqref{eq:3.7}, we obtain
$$
q(G)\leq q(G_3)<2n-2l<q(K_{bl-1}\vee(K_{n-(b+1)l+2}\cup(l-1)K_1)),
$$
which contradicts $q(G)\geq q(K_{bl-1}\vee(K_{n-(b+1)l+2}\cup(l-1)K_1))$. This completes the proof of Theorem 1.1. \hfill $\Box$

\section{The proof of Theorem 1.2}

\noindent{\it Proof of Theorem 1.2.} Suppose that $t_l(G)<\frac{1}{b}$ with $l\geq2$. Then we conclude $|S|<\frac{1}{b}c(G-S)$ for some nonempty subset $S\subseteq V(G)$. Notice that $c(G-S)\geq l$.
Thus, we obtain $c(G-S)\geq\max\{b|S|+1,l\}$. Let $|S|=s$ and $\omega=\max\{bs+1,l\}$. Obviously, $G$ is a spanning subgraph of $G_1=K_s\vee(K_{n_1}\cup K_{n_2}\cup\cdots\cup K_{n_{\omega}})$, where
$n_1\geq n_2\geq\cdots\geq n_{\omega}$ are positive integers with $\sum\limits_{i=1}^{\omega}n_i=n-s$. Write $G_2=K_s\vee(K_{n-s-\omega+1}\cup(\omega-1)K_1)$. Using Lemma 2.2 and Lemma 2.3, we infer
\begin{align}\label{eq:4.1}
q(G)\leq q(G_1)\leq q(G_2),
\end{align}
with equality occurring if and only if $G=G_1=G_2$.

\noindent{\bf Case 1.} $bs+1\geq l$.

Obviously, $\omega=bs+1$, $s\geq\lceil\frac{l-1}{b}\rceil$ and $G_2=K_s\vee(K_{n-bs-s}\cup bsK_1)$. Let $G_3=K_{\lceil\frac{l-1}{b}\rceil}\vee(K_{n-b\lceil\frac{l-1}{b}\rceil-\lceil\frac{l-1}{b}\rceil}\cup b\lceil\frac{l-1}{b}\rceil K_1)$.
If $s=\lceil\frac{l-1}{b}\rceil$, then $G_2=G_3$, and so $q(G_2)=q(G_3)$. Next, we discuss $s\geq\lceil\frac{l-1}{b}\rceil+1$.

For the partition $V(G_2)=V(K_s)\cup V(K_{n-bs-s})\cup V(bsK_1)$, the quotient matrix of $Q(G_2)$ is
\begin{align*}
B_1=\left(
  \begin{array}{ccc}
  n+s-2 & n-bs-s & bs\\
  s & 2n-2bs-s-2 & 0\\
  s & 0 & s\\
  \end{array}
\right).
\end{align*}
The characteristic polynomial of $B_1$ equals
\begin{align*}
f_{B_1}(x)=&x^{3}+(-3n+2bs-s+4)x^{2}\nonumber\\
&+(2n^{2}-2bsn+3sn-6n-4bs^{2}+4bs-4s+4)x\nonumber\\
&-2sn^{2}+4bs^{2}n+6sn-2b^{2}s^{3}-6bs^{2}-4s.
\end{align*}

Recall that $G_3=K_{\lceil\frac{l-1}{b}\rceil}\vee(K_{n-b\lceil\frac{l-1}{b}\rceil-\lceil\frac{l-1}{b}\rceil}\cup b\lceil\frac{l-1}{b}\rceil K_1)$. For the partition
$V(G_3)=V(K_{\lceil\frac{l-1}{b}\rceil})\cup V(K_{n-b\lceil\frac{l-1}{b}\rceil-\lceil\frac{l-1}{b}\rceil})\cup V(b\lceil\frac{l-1}{b}\rceil K_1)$, the quotient matrix of $Q(G_3)$ equals
\begin{align*}
B_2=\left(
  \begin{array}{ccc}
  n+\lceil\frac{l-1}{b}\rceil-2 & n-b\lceil\frac{l-1}{b}\rceil-\lceil\frac{l-1}{b}\rceil & b\lceil\frac{l-1}{b}\rceil\\
  \lceil\frac{l-1}{b}\rceil & 2n-2b\lceil\frac{l-1}{b}\rceil-\lceil\frac{l-1}{b}\rceil-2 & 0\\
  \lceil\frac{l-1}{b}\rceil & 0 & \lceil\frac{l-1}{b}\rceil\\
  \end{array}
\right),
\end{align*}
and so its characteristic polynomial is
\begin{align*}
f_{B_2}(x)=&x^{3}+\left(-3n+2b\left\lceil\frac{l-1}{b}\right\rceil-\left\lceil\frac{l-1}{b}\right\rceil+4\right)x^{2}\nonumber\\
&+\left(2n^{2}-2b\left\lceil\frac{l-1}{b}\right\rceil n+3\left\lceil\frac{l-1}{b}\right\rceil n-6n-4b\left\lceil\frac{l-1}{b}\right\rceil^{2}+4b\left\lceil\frac{l-1}{b}\right\rceil-4\left\lceil\frac{l-1}{b}\right\rceil+4\right)x\nonumber\\
&-2\left\lceil\frac{l-1}{b}\right\rceil n^{2}+4b\left\lceil\frac{l-1}{b}\right\rceil^{2}n+6\left\lceil\frac{l-1}{b}\right\rceil n-2b^{2}\left\lceil\frac{l-1}{b}\right\rceil^{3}-6b\left\lceil\frac{l-1}{b}\right\rceil^{2}-4\left\lceil\frac{l-1}{b}\right\rceil.
\end{align*}
By a simple calculation, we possess
\begin{align}\label{eq:4.2}
f_{B_1}(x)-f_{B_2}(x)=\left(s-\left\lceil\frac{l-1}{b}\right\rceil\right)\varphi(x),
\end{align}
where $\varphi(x)=(2b-1)x^{2}-(2bn-3n+4bs+4b\lceil\frac{l-1}{b}\rceil-4b+4)x-2n^{2}+4bsn+4bn\lceil\frac{l-1}{b}\rceil+6n-2b^{2}s^{2}-2b^{2}s\lceil\frac{l-1}{b}\rceil-2b^{2}\lceil\frac{l-1}{b}\rceil^{2}-6bs
-6b\lceil\frac{l-1}{b}\rceil-4$. The symmetry axis of $\varphi(x)$ is
\begin{align*}
x=&\frac{2bn-3n+4bs+4b\lceil\frac{l-1}{b}\rceil-4b+4}{2(2b-1)}\\
=&2n-2b\left\lceil\frac{l-1}{b}\right\rceil-2\\
&-\frac{1}{4b-2}\left(4n-4bs-4\left\lceil\frac{l-1}{b}\right\rceil-4+(6b-5)n-8b^{2}\left\lceil\frac{l-1}{b}\right\rceil+4b\left\lceil\frac{l-1}{b}\right\rceil-4b+4\right)\\
<&2n-2b\left\lceil\frac{l-1}{b}\right\rceil-2
\end{align*}
by $s\geq\lceil\frac{l-1}{b}\rceil+1$, $n\geq bs+s+1$ and $n\geq6b\left\lceil\frac{l-1}{b}\right\rceil>\frac{8b^{2}\lceil\frac{l-1}{b}\rceil-4b\lceil\frac{l-1}{b}\rceil+4b-4}{6b-5}$. This implies that $\varphi(x)$ is increasing for
$x\geq2n-2b\lceil\frac{l-1}{b}\rceil-2$. Thus, we have
\begin{align}\label{eq:4.3}
\varphi(x)\geq&\varphi\left(2n-2b\left\lceil\frac{l-1}{b}\right\rceil-2\right)\nonumber\\
=&(2b-1)\left(2n-2b\left\lceil\frac{l-1}{b}\right\rceil-2\right)^{2}\nonumber\\
&-\left(2bn-3n+4bs+4b\left\lceil\frac{l-1}{b}\right\rceil-4b+4\right)\left(2n-2b\left\lceil\frac{l-1}{b}\right\rceil-2\right)\nonumber\\
&-2n^{2}+4bsn+4bn\left\lceil\frac{l-1}{b}\right\rceil+6n-2b^{2}s^{2}-2b^{2}s\left\lceil\frac{l-1}{b}\right\rceil\nonumber\\
&-2b^{2}\left\lceil\frac{l-1}{b}\right\rceil^{2}-6bs-6b\left\lceil\frac{l-1}{b}\right\rceil-4\nonumber\\
=&-2b^{2}s^{2}+\left(-4bn+6b^{2}\left\lceil\frac{l-1}{b}\right\rceil+2b\right)s+4bn^{2}-12b^{2}n\left\lceil\frac{l-1}{b}\right\rceil\nonumber\\
&-2bn\left\lceil\frac{l-1}{b}\right\rceil-4bn+8b^{3}\left\lceil\frac{l-1}{b}\right\rceil^{2}+2b^{2}\left\lceil\frac{l-1}{b}\right\rceil^{2}+8b^{2}\left\lceil\frac{l-1}{b}\right\rceil+2b\left\lceil\frac{l-1}{b}\right\rceil\nonumber\\
\geq&-2b^{2}\left(\frac{n-1}{b+1}\right)^{2}+\left(-4bn+6b^{2}\left\lceil\frac{l-1}{b}\right\rceil+2b\right)\left(\frac{n-1}{b+1}\right)+4bn^{2}\nonumber\\
&-12b^{2}n\left\lceil\frac{l-1}{b}\right\rceil-2bn\left\lceil\frac{l-1}{b}\right\rceil-4bn+8b^{3}\left\lceil\frac{l-1}{b}\right\rceil^{2}+2b^{2}\left\lceil\frac{l-1}{b}\right\rceil^{2}\nonumber\\
&+8b^{2}\left\lceil\frac{l-1}{b}\right\rceil+2b\left\lceil\frac{l-1}{b}\right\rceil \ \ \ \ \ \ \ \ \left(\mbox{since} \ s\leq\frac{n-1}{b+1}\right)\nonumber\\
=&\frac{(4b^{3}+2b^{2})n^{2}}{(b+1)^{2}}-\frac{((12b^{4}+20b^{3}+10b^{2}+2b)\lceil\frac{l-1}{b}\rceil+4b^{3}-2b^{2}-2b)n}{(b+1)^{2}}\nonumber\\
&+\frac{(8b^{5}+18b^{4}+12b^{3}+2b^{2})\lceil\frac{l-1}{b}\rceil^{2}+(8b^{4}+12b^{3}+6b^{2}+2b)\lceil\frac{l-1}{b}\rceil-4b^{2}-2b}{(b+1)^{2}}\nonumber\\
\geq&\frac{(80b^{5}-30b^{4}-48b^{3}-10b^{2})\lceil\frac{l-1}{b}\rceil^{2}}{(b+1)^{2}}\nonumber\\
&+\frac{(-16b^{4}+24b^{3}+6b^{2}+2b)\lceil\frac{l-1}{b}\rceil-4b^{2}-2b}{(b+1)^{2}} \ \ \ \ \ \ \ \ \left(\mbox{since} \ n\geq6b\left\lceil\frac{n-1}{b+1}\right\rceil\right)\nonumber\\
\geq&\frac{(80b^{5}-30b^{4}-48b^{3}-10b^{2})\lceil\frac{l-1}{b}\rceil}{(b+1)^{2}}\nonumber\\
&+\frac{(-16b^{4}+24b^{3}+6b^{2}+2b)\lceil\frac{l-1}{b}\rceil-4b^{2}-2b}{(b+1)^{2}} \ \ \ \ \ \ \ \ \left(\mbox{since} \ b\geq2 \ \mbox{and} \ \left\lceil\frac{n-1}{b+1}\right\rceil\geq1\right)\nonumber\\
=&\frac{(80b^{5}-48b^{4}-24b^{3}-4b^{2}+2b)\lceil\frac{l-1}{b}\rceil-4b^{2}-2b}{(b+1)^{2}}\nonumber\\
\geq&\frac{80b^{5}-48b^{4}-24b^{3}-8b^{2}}{(b+1)^{2}} \ \ \ \ \ \ \ \ \left(\mbox{since} \ b\geq2 \ \mbox{and} \ \left\lceil\frac{n-1}{b+1}\right\rceil\geq1\right)\nonumber\\
>&0 \ \ \ \ \ \ \ \ (\mbox{since} \ b\geq2).
\end{align}
In terms of \eqref{eq:4.2}, \eqref{eq:4.3} and $s\geq\lceil\frac{l-1}{b}\rceil+1$, we conclude $f_{B_1}(x)>f_{B_2}(x)$ for $x\geq2n-2b\lceil\frac{l-1}{b}\rceil-2$. Note that
$G_3=K_{\lceil\frac{l-1}{b}\rceil}\vee(K_{n-b\lceil\frac{l-1}{b}\rceil-\lceil\frac{l-1}{b}\rceil}\cup b\lceil\frac{l-1}{b}\rceil K_1)$ contains $K_{n-b\lceil\frac{l-1}{b}\rceil}$ as a proper subgraph. Then it
follows from Lemma 2.2 that $q(G_3)>q(K_{n-b\lceil\frac{l-1}{b}\rceil})=2n-2b\lceil\frac{l-1}{b}\rceil-2$, and so $q(G_3)>q(G_2)$ by Lemma 2.1.

From the above discussion, we possess
\begin{align}\label{eq:4.4}
q(G_2)\leq q(G_3),
\end{align}
with equality if and only if $G_2=G_3$.

Let $G_*=K_{\lfloor\frac{l-1}{b}\rfloor}\vee(K_{n-\lfloor\frac{l-1}{b}\rfloor-l+1}\cup(l-1)K_1)$. Notice that if $b|(l-1)$, then $G_3=G_*$, and so $q(G_3)=q(G_*)$. Let $b\nmid(l-1)$ and
$G_3'=K_{\lceil\frac{l-1}{b}\rceil-1}\vee(K_{n-b\lceil\frac{l-1}{b}\rceil-\lceil\frac{l-1}{b}\rceil+2}\cup(b\lceil\frac{l-1}{b}\rceil-1)K_1)$. For the partition $V(G_3')=V(K_{\lceil\frac{l-1}{b}\rceil-1})
\cup V(K_{n-b\lceil\frac{l-1}{b}\rceil-\lceil\frac{l-1}{b}\rceil+2})\cup V((b\lceil\frac{l-1}{b}\rceil-1)K_1)$, the quotient matrix of $Q(G_3')$ equals
\begin{align*}
B_2'=\left(
  \begin{array}{ccc}
  n+\lceil\frac{l-1}{b}\rceil-3 & n-b\lceil\frac{l-1}{b}\rceil-\lceil\frac{l-1}{b}\rceil+2 & b\lceil\frac{l-1}{b}\rceil-1\\
  \lceil\frac{l-1}{b}\rceil-1 & 2n-2b\lceil\frac{l-1}{b}\rceil-\lceil\frac{l-1}{b}\rceil+1 & 0\\
  \lceil\frac{l-1}{b}\rceil-1 & 0 & \lceil\frac{l-1}{b}\rceil-1\\
  \end{array}
\right).
\end{align*}
The characteristic polynomial of $B_2'$ is
\begin{align*}
f_{B_2'}(x)=&x^{3}+\left(-3n+2b\left\lceil\frac{l-1}{b}\right\rceil-\left\lceil\frac{l-1}{b}\right\rceil+3\right)x^{2}\nonumber\\
&+\left(2n^{2}-2b\left\lceil\frac{l-1}{b}\right\rceil n+3\left\lceil\frac{l-1}{b}\right\rceil n-7n-4b\left\lceil\frac{l-1}{b}\right\rceil^{2}+8b\left\lceil\frac{l-1}{b}\right\rceil\right)x\nonumber\\
&-2\left\lceil\frac{l-1}{b}\right\rceil n^{2}+2n^{2}+4b\left\lceil\frac{l-1}{b}\right\rceil^{2}n-4b\left\lceil\frac{l-1}{b}\right\rceil n+2\left\lceil\frac{l-1}{b}\right\rceil n-2n\nonumber\\
&-2b^{2}\left\lceil\frac{l-1}{b}\right\rceil^{3}+2b^{2}\left\lceil\frac{l-1}{b}\right\rceil^{2}-2b\left\lceil\frac{l-1}{b}\right\rceil^{2}+2b\left\lceil\frac{l-1}{b}\right\rceil.
\end{align*}
For $x\geq2n-2b\lceil\frac{l-1}{b}\rceil$, we deduce
\begin{align}\label{eq:4.5}
f_{B_2}(x)-f_{B_2'}(x)=&x^{2}+\left(n-4b\left\lceil\frac{l-1}{b}\right\rceil-4\left\lceil\frac{l-1}{b}\right\rceil+4\right)x\nonumber\\
&-2n^{2}+4b\left\lceil\frac{l-1}{b}\right\rceil n+4\left\lceil\frac{l-1}{b}\right\rceil n+2n-2b^{2}\left\lceil\frac{l-1}{b}\right\rceil^{2}\nonumber\\
&-4b\left\lceil\frac{l-1}{b}\right\rceil^{2}-2b\left\lceil\frac{l-1}{b}\right\rceil-4\left\lceil\frac{l-1}{b}\right\rceil\nonumber\\
\geq&\left(2n-2b\left\lceil\frac{l-1}{b}\right\rceil\right)^{2}+\left(n-4b\left\lceil\frac{l-1}{b}\right\rceil-4\left\lceil\frac{l-1}{b}\right\rceil+4\right)\left(2n-2b\left\lceil\frac{l-1}{b}\right\rceil\right)\nonumber\\
&-2n^{2}+4b\left\lceil\frac{l-1}{b}\right\rceil n+4\left\lceil\frac{l-1}{b}\right\rceil n+2n-2b^{2}\left\lceil\frac{l-1}{b}\right\rceil^{2}\nonumber\\
&-4b\left\lceil\frac{l-1}{b}\right\rceil^{2}-2b\left\lceil\frac{l-1}{b}\right\rceil-4\left\lceil\frac{l-1}{b}\right\rceil \ \ \ \ \ \ \ \left(\mbox{since} \ x\geq2n-2b\left\lceil\frac{l-1}{b}\right\rceil\right)\nonumber\\
=&4n^{2}-\left(14b\left\lceil\frac{l-1}{b}\right\rceil+4\left\lceil\frac{l-1}{b}\right\rceil-10\right)n+10b^{2}\left\lceil\frac{l-1}{b}\right\rceil^{2}\nonumber\\
&+4b\left\lceil\frac{l-1}{b}\right\rceil^{2}-10b\left\lceil\frac{l-1}{b}\right\rceil-4\left\lceil\frac{l-1}{b}\right\rceil\nonumber\\
\geq&4\left(6b\left\lceil\frac{l-1}{b}\right\rceil\right)^{2}-\left(14b\left\lceil\frac{l-1}{b}\right\rceil+4\left\lceil\frac{l-1}{b}\right\rceil-10\right)\left(6b\left\lceil\frac{l-1}{b}\right\rceil\right)
+10b^{2}\left\lceil\frac{l-1}{b}\right\rceil^{2}\nonumber\\
&+4b\left\lceil\frac{l-1}{b}\right\rceil^{2}-10b\left\lceil\frac{l-1}{b}\right\rceil-4\left\lceil\frac{l-1}{b}\right\rceil \ \ \ \ \ \ \ \left(\mbox{since} \ n\geq6b\left\lceil\frac{l-1}{b}\right\rceil\right)\nonumber\\
=&(70b^{2}-20b)\left\lceil\frac{l-1}{b}\right\rceil^{2}+(50b-4)\left\lceil\frac{l-1}{b}\right\rceil\nonumber\\
>&0 \ \ \ \ \ \ \ \left(\mbox{since} \ b\geq2 \ \mbox{and} \ \left\lceil\frac{l-1}{b}\right\rceil\geq1\right).
\end{align}
Notice that $G_3'$ contains $K_{n-b\lceil\frac{l-1}{b}\rceil+1}$ as a proper subgraph. According to Lemma 2.2, we obtain $q(G_3')>q(K_{n-b\lceil\frac{l-1}{b}\rceil+1})=2n-2b\lceil\frac{l-1}{b}\rceil$,
and so $q(G_3)<q(G_3')$ by \eqref{eq:4.5} and Lemma 2.1. Notice that $\lfloor\frac{l-1}{b}\rfloor=\lceil\frac{l-1}{b}\rceil-1$ and $b\lceil\frac{l-1}{b}\rceil-1\geq l-1$. Then we see that $G_3'$ is a spanning
subgraph of $G_*$, and so $q(G_3')\leq q(G_*)$. From the above discussion, we have $q(G_3)\leq q(G_*)$, with equality if and only if $G_3=G_*$.Combining this with \eqref{eq:4.1} and \eqref{eq:4.4}, we infer
$$
q(G)\leq q(G_*)=q(K_{\lfloor\frac{l-1}{b}\rfloor}\vee(K_{n-\lfloor\frac{l-1}{b}\rfloor-l+1}\cup(l-1)K_1)),
$$
where the equality occurs if and only if $G=K_{\lfloor\frac{l-1}{b}\rfloor}\vee(K_{n-\lfloor\frac{l-1}{b}\rfloor-l+1}\cup(l-1)K_1)$, which contradicts the condition of the theorem because
$K_{\lfloor\frac{l-1}{b}\rfloor}\vee(K_{n-\lfloor\frac{l-1}{b}\rfloor-l+1}\cup(l-1)K_1)$ is not $(\frac{1}{b},l)$-tough.

\noindent{\bf Case 2.} $bs+1<l$.

Obviously, $\omega=l$ and $s<\frac{l-1}{b}$. Recall that $G_2=K_s\vee(K_{n-s-\omega+1}\cup(\omega-1)K_1)$. Notice that $G_2$ is a spanning subgraph of $K_{\lfloor\frac{l-1}{b}\rfloor}\vee(K_{n-\lfloor\frac{l-1}{b}\rfloor-l+1}\cup(l-1)K_1)$.
Combining this with \eqref{eq:4.1} and Lemma 2.2, we conclude
$$
q(G)\leq q(G_2)\leq q(K_{\lfloor\frac{l-1}{b}\rfloor}\vee(K_{n-\lfloor\frac{l-1}{b}\rfloor-l+1}\cup(l-1)K_1)),
$$
with equality if and only if $G=K_{\lfloor\frac{l-1}{b}\rfloor}\vee(K_{n-\lfloor\frac{l-1}{b}\rfloor-l+1}\cup(l-1)K_1)$, which contradicts the condition of the theorem because
$K_{\lfloor\frac{l-1}{b}\rfloor}\vee(K_{n-\lfloor\frac{l-1}{b}\rfloor-l+1}\cup(l-1)K_1)$ is not $(\frac{1}{b},l)$-tough. This completes the proof of Theorem 1.2. \hfill $\Box$

\section*{Declaration of competing interest}

The authors declare that they have no known competing financial interests or personal relationships that could have appeared to influence the work reported in this paper.

\section*{Data availability}

No data was used for the research described in the article.

\section*{Acknowledgments}

This work was supported by the Natural Science Foundation of Jiangsu Province (Grant No. BK20241949). Project ZR2023MA078 supported by Shandong Provincial Natural Science Foundation.

\end{document}